\documentclass[a4paper,11pt,twoside]{article}

\RequirePackage[a4paper,head=1cm]{geometry}
\usepackage[utf8x]{inputenc}
\usepackage{amsmath}
\usepackage{amssymb}
\usepackage{hyperref}
\usepackage{dsfont}
\newtheorem{Definition}{Definition}[section]
\newtheorem{Theoreme}{Theorem}
\newtheorem{Lemme}{Lemma}[section]

\newtheorem{Proposition}{Proposition}[section]
\newtheorem{Remarque}{\bf Remark}

\title{\bf A counterexample for\\ Improved Sobolev Inequalities\\ over the $2$-adic group} 
\author{Diego Chamorro } 

\begin{document} 
\maketitle 
\begin{center}
\begin{minipage}[l]{150mm}
\begin{scriptsize}\abstract{On the framework of the $2$-adic group $\mathbb{Z}_2$, we study a Sobolev-like inequality where we estimate the $L^2$ norm by a geometric mean of the $BV$ norm and the $\dot{B}^{-1,\infty}_{\infty}$ norm. We first show, using the special topological properties of the $p$-adic groups, that the set of functions of bounded variations $BV$ can be identified to the Besov space $\dot{B}^{1,\infty}_{1}$. This identification lead us to the construction of a counterexample to the improved Sobolev inequality.\\[3mm]
}\end{scriptsize}
\end{minipage}
\end{center}
\section{Introduction}

The general improved Sobolev inequalities were initially introduced by P. Gérard, Y. Meyer and F. Oru in \cite{GMO}. For a function $f$ such that $f\in \dot{W}^{s_1,p}(\mathbb{R}^n)$ and $f\in \dot{B}^{-\beta,\infty}_{\infty}(\mathbb{R}^n)$, these inequalities read as follows: 
\begin{equation}\label{GeneralISP}
\|f\|_{\dot{W}^{s,q}}\leq C \|f\|_{\dot{W}^{s_1,p}}^\theta\|f\|_{\dot{B}^{-\beta,\infty}_{\infty}}^{1-\theta}
\end{equation}
where $1<p<q<+\infty$, $\theta=p/q$, $s=\theta s_1 -(1-\theta)\beta$ and $-\beta<s<s_1$. The method used for proving these estimates relies on the Littlewood-Paley decomposition and on a dyadic bloc manipulation and this explains the fact that the value $p=1$ is forbidden here.\\

In order to study the case $p=1$, it is necessary to develop other techniques. The case when $p=1$, $s=0$ and $s_1=1$ was treated by M. Ledoux in \cite{Ledoux} using a special cut-off function; while the case  $s_1=1$ and $p=1$ was studied by  A. Cohen, W. Dahmen, I. Daubechies \& R. De Vore in \cite{Cohen2}. In this last article, the authors give a BV-norm weak estimation using wavelet coefficients and isoperimetric inequalities and obtained, for a function $f$ such that $f\in BV(\mathbb{R}^n)$ and $f\in \dot{B}^{-\beta,\infty}_{\infty}(\mathbb{R}^n)$, the estimation below:
\begin{equation}\label{ISI2}
\|f\|_{\dot{W}^{s,q}}\leq C \| f\|_{BV}^{1/q}\|f\|_{\dot{B}^{-\beta,\infty}_{\infty}}^{1-1/q}
\end{equation}
where $1<q\leq 2$, $0\leq s<1/q$ and $\beta=(1-sq)/(q-1)$. \\

In a previous work (see \cite{Chame}, \cite{Chame1}), we studied the possible generalizations of inequalities of type (\ref{GeneralISP}) and (\ref{ISI2}) to other frameworks than $\mathbb{R}^n$. In particular, we worked over stratified Lie groups and over polynomial volume growth Lie groups and we obtained some new weak-type estimates. \\

The aim of this paper is to study inequalities of type (\ref{GeneralISP}) and (\ref{ISI2}) in the setting of the $2$-adic group $\mathbb{Z}_2$.
The main reason for working in the framework of $\mathbb{Z}_2$ is that this group is completely different from $\mathbb{R}^n$ and from stratified or polynomial Lie groups. Indeed, since the $2$-adic group is totally discontinuous, it is not absolutely trivial to give a definition for smoothness measuring spaces. Thus, the first step to do, in order to study these Sobolev-like inequalities, is to give an adapted characterization of such functional spaces. This will be achieved using the Littlewood-Paley approach and, once this task is done, we will immediatly prove -following the classical path exposed in \cite{GMO}- the inequalities  (\ref{GeneralISP}) in the setting of the $2$-adic group $\mathbb{Z}_2$.\\

For the estimate (\ref{ISI2}), we introduce the $BV$ space in the following manner:  we will say that $f\in BV(\mathbb{Z}_{2})$ if there exists a constant $C>0$ such that
$$\int_{\mathbb{Z}_{2}}|f(x+y)-f(x)|dx \leq C|y|_{2}\quad  (\forall y \in \mathbb{Z}_{2}).$$
As a surprising fact, we obtain the 
\begin{Theoreme}\label{IdentityBVBesov}
We have the following relationship between the space of functions of bounded variation $BV(\mathbb{Z}_2)$ and the Besov space $\dot{B}^{1,\infty}_1(\mathbb{Z}_2)$: 
$$BV(\mathbb{Z}_2)\simeq \dot{B}^{1,\infty}_1(\mathbb{Z}_2)$$
\end{Theoreme}

Of course, this identification is false in $\mathbb{R}^n$ and it is this special relationship in $\mathbb{Z}_2$ that give us our principal theorem which is the $2$-adic counterpart of the inequality (\ref{ISI2}):

\begin{Theoreme}\label{MainTheorem}
The following inequality is false in $\mathbb{Z}_{2}$. There is not an universal constant $C>0$ such that we have
$$\|f\|_{L^2}^{2}\leq C\|f\|_{BV}\|f\|_{\dot{B}^{-1,\infty}_{\infty}}$$
for all $f\in BV \cap \dot{B}^{-1,\infty}_{\infty}(\mathbb{Z}_{2})$.
\end{Theoreme}
 
This striking fact says that the improved Sobolev inequalities of type (\ref{ISI2}) depend on the group's structure and that they are no longer true for the $2$-adic group $\mathbb{Z}_{2}$.\\

The plan of the article is the following: in section \ref{IntroP} we recall some well known properties about $p$-adic groups, in \ref{ESPFUNC} we define Sobolev and Besov spaces, in \ref{BVadique0} we prove theorem \ref{IdentityBVBesov} and, finally, we prove the theorem \ref{MainTheorem} in section \ref{ISPadique00}.

\section{$p$-adic groups}\label{IntroP}

We write $a|b$ when $a$ divide $b$ or, equivalently, when $b$ is a multiple of $a$. Let $p$ be any prime number, for $0\neq x\in \mathbb{Z}$, we define the $p$-adic valuation of $x$ by $\gamma(x)=\max\{r: p^{r}|x\}\geq 0$ and, for any rational number $x=\frac{a}{b}\in \mathbb{Q}$, we write $\gamma(x)=\gamma(a)-\gamma(b)$. Furthermore if $x=0$, we agree to write $\gamma(0)=+\infty$.\\ 

Let $x \in \mathbb{Q}$ and $p$ be any prime number, with the $p$-adic valuation of $x$ we can construct a norm by writing
\begin{equation}\label{Normeadique}
|x|_{p}=\left\lbrace \begin{array}{ll}
p^{-\gamma} &\mbox{if}\quad x\neq 0\\[5mm]
p^{-\infty}=0 &\mbox{if}\quad x= 0.
\end{array}\right.
\end{equation}
This expression satisfy the following properties
\begin{enumerate}
\item[a)]$|x|_{p}\geq 0$, and $|x|_{p}=0 \iff x=0$;
\item[b)] $|xy|_{p}=|x|_{p}|y|_{p}$;
\item[c)] $|x+y|_{p}\leq \max\{|x|_{p},|y|_{p}\}$, with equality when $|x|_{p}\neq |y|_{p}$.
\end{enumerate}
When a norm satisfy $c)$ it is called a non-Archimedean norm and an interesting fact is that over $\mathbb{Q}$ \textit{all} the possible norms are equivalent to $|\cdot|_{p}$ for some $p$: this is the so-called Ostrowski theorem, see \cite{Amice} for a proof. 
\begin{Definition}
Let $p$ be a any prime number. We define the field of $p$-adic numbers $\mathbb{Q}_{p}$ as the completion of $\mathbb{Q}$ when using the norm $|\cdot|_{p}$.
\end{Definition}
We present in the following lines the algebraic structure of the set $\mathbb{Q}_{p}$. Every $p$-adic number $x\neq 0$ can be represented in a unique manner by the formula
\begin{equation}\label{Hensel}
x=p^{\gamma}(x_{0}+x_{1}p+x_{2}p^{2}+...),
\end{equation}
where $\gamma=\gamma(x)$ is the $p$-adic valuation of $x$ and $x_{j}$ are integers such that $x_{0}>0$ and $0\leq x_{j}\leq p-1$ for $j=1,2,...$. Remark that this canonical representation implies the identity $|x|_{p}=p^{-\gamma}$.\\

Let $x,y \in \mathbb{Q}_{p}$, using the formula (\ref{Hensel}) we define the sum of $x$ and $y$ by 
$x+y=p^{\gamma (x+y)}(c_{0}+c_{1}p+c_{2}p^{2}+...)$  with $0\leq c_{j}\leq p-1$ and $c_{0}>0$, where $\gamma(x+y)$ and $c_{j}$ are the unique solution of the equation
$$p^{\gamma(x)}(x_{0}+x_{1}p+x_{2}p^{2}+...)+ p^{\gamma(y)}(y_{0}+y_{1}p+y_{2}p^{2}+...)=p^{\gamma(x+y)}(c_{0}+c_{1}p+c_{2}p^{2}+...).$$
Furthermore, for $a,x \in \mathbb{Q}_{p}$, the equation $a+x=0$ has a unique solution in $\mathbb{Q}_{p}$ given by $x=-a$. In the same way, the equation $ax=1$ has a unique solution in $\mathbb{Q}_{p}$: $x=1/a$.\\

We take now a closer look at the topological structure of $\mathbb{Q}_{p}$. With the norm $|\cdot|_p$ we construct a distance over $\mathbb{Q}_p$ by writing
\begin{equation}\label{dist}
d(x,y)=|x-y|_{p} 
\end{equation} 
and we define the balls $B_{\gamma}(x)=\left\{y\in \mathbb{Q}_{p}: \; d(x,y)\leq p^{\gamma}\right\}$  with $\gamma \in \mathbb{Z}$. Remark that, from the properties of the $p$-adic valuation, this distance has the \textit{ultra-metric} property (\textit{i.e.} $d(x,y)\leq \max\{d(x,z),d(z,y)\}\leq |x|_{p}+|y|_{p}$).\\

We gather with the next proposition some important facts concerning the balls in $\mathbb{Q}_{p}$.
\begin{Proposition}
Let $\gamma$ be an integer, then we have
\begin{enumerate}
\item[1)] the ball $B_{\gamma}(x)$ is a open and a closed set for the distance (\ref{dist}).
\item[2)] every point of $B_{\gamma}(x)$ is its center.
\item[3)] $\mathbb{Q}_{p}$ endowed with this distance is a complete Hausdorff metric space.
\item[4)] $\mathbb{Q}_{p}$ is a locally compact set.
\item[5)] the $p$-adic group $\mathbb{Q}_{p}$ is a totally discontinuous space.
\end{enumerate}
\end{Proposition}

For a proof of this proposition and more details see the books \cite{Amice}, \cite{Koblitz} or \cite{VVZ}. 
\section{Functional spaces}\label{ESPFUNC}
In this article, we will work with the subset $\mathbb{Z}_{2}$ of $\mathbb{Q}_{2}$ which is defined by $\mathbb{Z}_{2}=\{x\in \mathbb{Q}_{2}:\; |x|_{2}\leq 1\}$, and we will focus on real-valued functions over $\mathbb{Z}_{2}$. Since $\mathbb{Z}_{2}$ is a locally compact commutative group, there exists a Haar measure $dx$ which is translation invariant $i.e.$: $d(x+a)=dx$, furthermore we have the identity $d(xa)=|a|_{2}dx$ for $a\in \mathbb{Z}_{2}^{*}$. We will normalize the measure $dx$ by setting
$$\int_{\{|x|_{2}\leq1\}}dx=1.$$
This measure is then unique and we will note $|E|$ the measure for any subset $E$ of $\mathbb{Z}_{2}$. Lebesgue spaces $L^{p}(\mathbb{Z}_{2})$ are thus defined in a natural way: $\|f\|_{L^p}=\left(\int_{\mathbb{Z}_{2}}|f(x)|^{p}dx\right)^{1/p}$ for $1\leq p< +\infty$, with the usual modifications when $p=+\infty$.\\

Let us now introduce the Littlewood-Paley decomposition in $\mathbb{Z}_{2}$. We note $\mathcal{F}_{j}$ the Boole algebra formed by the equivalence classes $E\subset \mathbb{Z}_{2}$ modulo the sub-group $2^{j}\mathbb{Z}_{2}$. Then, for any function $f\in L^{1}(\mathbb{Z}_{2})$, we call $S_{j}(f)$ the conditionnal expectation of $f$ with respect to $\mathcal{F}_{j}$:
$$S_{j}(f)(x)=\frac{1}{|B_{j}(x)|}\int_{B_{j}(x)}f(y)dy.$$
The dyadic blocks are thus defined by the formula $\Delta_{j}(f)=S_{j+1}(f)-S_{j}(f)$ and the Littlewood-Paley decomposition of a function $f:\mathbb{Z}_{2}\longrightarrow \mathbb{R}$ is given by
\begin{equation}\label{LPP}
f=S_{0}(f)+\sum_{j=0}^{+\infty}\Delta_{j}(f)\qquad \mbox{where } S_{0}(f)=\int_{\mathbb{Z}_{2}}f(x)dx.
\end{equation}

We will need in the sequel some very special sets noted $Q_{j,k}$. Here is the definition and some properties:
\begin{Proposition}\label{pala2bis}
Let $j\in \mathbb{N}$ and $k=\{0,1,...,2^{j}-1\}$. Define the subset $Q_{j,k}$ of $\mathbb{Z}_2$ by
\begin{equation}\label{QJKset}
Q_{j,k}=\left\{k+2^{j}\mathbb{Z}_{2}\right\}.
\end{equation}
Then
\begin{enumerate}
\item[1)] We have the identity $\mathcal{F}_{j}=\underset{0\leq k<2^{j}}{\bigcup} Q_{j,k}$,
\item[2)] For $k=\{0,1,...,2^{j}-1\}$ the sets $Q_{j,k}$ are mutually disjoint, 
\item[3)] $|Q_{j,k}|=2^{-j}$ for all $k$,
\item[4)] the $2$-adic valuation is constant over $Q_{j,k}$.
\end{enumerate}
\end{Proposition}

The verifications are easy and left to the reader.\\

With the Littlewood-Paley decomposition given in (\ref{LPP}), we obtain the following equivalence for the Lebesgue spaces $L^p(\mathbb{Z}_2)$ with $1<p<+\infty$:
$$\|f\|_{L^p}\simeq \|S_0(f)\|_{L^p}+\left\|\bigg(\sum_{j\in \mathbb{N}}|\Delta_{j}f|^{2}\bigg)^{1/2} \right\|_{L^p}.$$
See the book \cite{Stein0}, chapter IV, for a general proof.\\

Let us turn now to smoothness measuring spaces. As said in the introduction, it is not absolutely trivial to define Sobolev and Besov spaces over $\mathbb{Z}_2$ since we are working in a totally discontinuous setting. Here is an example of this situation with the Sobolev space $W^{1,2}$: one could try to define the quantity $|\nabla f|$ by the formula
\begin{equation*}
|\nabla f|=\underset{\delta \to 0}{\lim}\;\underset{d(x,y)<\delta}{\sup} \frac{|f(x)-f(y)|}{d(x,y)}
\end{equation*}
and define the Sobolev space $W^{1,2}(\mathbb{Z}_2)$ by the norm
\begin{equation}\label{padique1}
\|f\|_{\ast}=\|f\|_{L^2}+\left(\int_{\mathbb{Z}_{2}}|\nabla f|^{2}dx\right)^{1/2}.
\end{equation}
Now, using the Littlewood-Paley decomposition we can also write
\begin{equation*}
\|f\|_{\ast \ast}=\|S_{0}f\|_{L^2}+\left\|\left(\sum_{j\in \mathbb{N}}2^{2j}|\Delta_{j}f|^{2}\right)^{1/2}\right\|_{2}.
\end{equation*}
However, the quantities $\|\cdot\|_{\ast}$ and $\|\cdot\|_{\ast\ast}$ are not equivalent: in the case of (\ref{padique1}) consider a function $f=c_{k}$ constant over each $Q_{j,k}=\{k+2^{j}\mathbb{Z}_{2}\}$ for some fixed $j$. Then  we have $|\nabla f|\equiv 0$ and for these functions the norm $\|\cdot\|_{\ast}$ would be equal to the $L^{2}$ norm.\\ 

This is the reason why we will use in this article the Littlewood-Paley approach to characterize Sobolev spaces:
\begin{equation}\label{Sobolevadique}
\|f\|_{W^{s,p}}\simeq\|S_{0}f\|_{L^p}+\left\|\bigg(\sum_{j\in \mathbb{N}}2^{2js}|\Delta_{j}f|^{2}\bigg)^{1/2} \right\|_{L^p}.
\end{equation}
with $1<p<+\infty$ and $s>0$. For Besov spaces we will define them by the norm
\begin{equation}\label{besov1}
\|f\|_{B^{s,q}_{p}}\simeq\|S_{0}f\|_{L^p}+\left(\sum_{j\in \mathbb{N}}2^{jsq}\|\Delta_{j}f\|^{q}_{L^p}\right)^{1/q}
\end{equation}
where $s\in \mathbb{R}$, $1\leq p,q<+\infty$ with the necessary modifications when $p,q=+\infty$.\\
\begin{Remarque}
\emph{For homogeneous functional spaces $\dot{W}^{s,p}$ and $\dot{B}^{s,q}_{p}$, we drop out the term $\|S_{0}f\|_{L^p}$ in (\ref{Sobolevadique}) and (\ref{besov1}).}
\end{Remarque}

Let us give some simple examples of function belonging to these functional spaces. 

\begin{enumerate}
\item[1)] The function $f(x)=\log_{2}|x|_{2}$ is in $\dot{B}^{1,\infty}_{1}(\mathbb{Z}_{2})$. First note that $|x|_{2}=2^{-\gamma(x)}$ and thus $f(x)=-\gamma(x)$. Recall (cf. proposition \ref{pala2bis}) that over each set $Q_{j,k}$, the quantity $\gamma(x)$ is constant, so the dyadic bloc $\Delta_{j}f$ is given by
$$\Delta_{j}f(x)=\left\lbrace\begin{array}{l}
-1 \qquad\mbox{over}\quad  Q_{j+1,0}\\[5mm]
0 \qquad\mbox{elsewhere}.
\end{array}\right.$$
Hence, taking the $L^{1}$ norm, we have $\|\Delta_{j}f\|_{L^1}= \frac{1}{2} 2^{-j}$ and then $f\in \dot{B}^{1,\infty}_{1}(\mathbb{Z}_{2})$.\\

\item[2)] Set $h(x)=1/|x|_{2}$, we have $h\in \dot{B}^{-1,\infty}_{\infty}$. For this, we must verify $\underset{j\geq 0}{\sup}\, 2^{-j}\|\Delta_{j}h\|_{L^\infty}<+\infty$. By definition we obtain $h(x)=2^{\gamma(x)}$ and then
$$\Delta_{j}h(x)=\left\lbrace\begin{array}{l}
2^{j} \qquad\mbox{over}\quad  Q_{j+1,0}\\[5mm]
0 \qquad\mbox{elsewhere}.
\end{array}\right.$$
We finally obtain $\|\Delta_{j}h\|_{L^\infty}=2^{j}$ and hence $2^{-j}\|\Delta_{j}h\|_{L^\infty}=1$ for all $j$, so we write $h\in \dot{B}^{-1,\infty}_{\infty}$.\\[5mm]
\end{enumerate}

With the Littlewood-Paley characterisation of Sobolev spaces and Besov spaces given in (\ref{Sobolevadique}) and (\ref{besov1}) we have the following theorem:

\begin{Theoreme} In the framework of the $2$-adic group $\mathbb{Z}_2$ we have, for a function $f$ such that $f\in \dot{W}^{s_1,p}(\mathbb{Z}_2)$ and $f\in \dot{B}^{-\beta,\infty}_{\infty}(\mathbb{Z}_2)$, the inequality
$$\|f\|_{\dot{W}^{s,q}}\leq C \|f\|_{\dot{W}^{s_1,p}}^\theta\|f\|_{\dot{B}^{-\beta,\infty}_{\infty}}^{1-\theta}$$
where $1<p<q<+\infty$, $\theta=p/q$, $s=\theta s_1 -(1-\theta)\beta$ and $-\beta<s<s_1$.
\end{Theoreme}
\textbf{\textit{Proof.}} We start with an interpolation result: let $(a_j)_{j\in \mathbb{N}}$ be a sequence, let $s=\theta s_1 -(1-\theta)\beta$ with $\theta=p/q$, then we have for $r,r_1,r_2\in [1,+\infty]$ the estimate
$$\|2^{js}a_j\|_{\ell^r}\leq C\|2^{js_1}a_j\|_{\ell^{r_1}}^\theta\|2^{-j\beta}a_j\|_{\ell^{r_2}}^{1-\theta}$$
See \cite{Bergh} for a proof. Apply this estimate to the dyadic blocks $\Delta_j f$ to obtain
$$\left(\sum_{j\in \mathbb{Z}}2^{2js}|\Delta_j f(x)|^2\right)^{1/2}\leq C \left(\sum_{j\in \mathbb{Z}}2^{2js_1}|\Delta_j f(x)|^2\right)^{\theta/2}\left(\underset{j\in \mathbb{Z}}{\sup}\; 2^{-j\beta}|\Delta_j f(x)| \right)^{1-\theta}$$
To finish, compute the $L^q$ norm of the preceding quantities. 

\begin{flushright}{$\blacksquare$}\end{flushright} 
\section{The $BV(\mathbb{Z}_{2})$ space and the proof of theorem \ref{IdentityBVBesov}}\label{BVadique0}

We study in this section the space of functions of bounded variation $BV$ and we will prove some surprising facts in the framework of $2$-adic group $\mathbb{Z}_{2}$. Let us start recalling the definition of this space:
\begin{Definition}\label{BVadique}
If $f$ is a real-valued function over $\mathbb{Z}_{2}$, we will say that $f\in BV(\mathbb{Z}_{2})$ if there exists a constant $C>0$ such that
\begin{equation}\label{BVadique1}
\int_{\mathbb{Z}_{2}}|f(x+y)-f(x)|dx \leq C|y|_{2},\quad  (\forall y \in \mathbb{Z}_{2}).
\end{equation}
\end{Definition}

We prove now the theorem \ref{IdentityBVBesov} which asserts that in $\mathbb{Z}_{2}$, the $BV$ space can be identified to the Besov space $\dot{B}^{1,\infty}_{1}$. For this, we will use two steps given by the propositions \ref{Dosadico} and \ref{Unoadico} below.

\begin{Proposition}\label{Dosadico}
If $f$ is a real-valued function over $\mathbb{Z}_{2}$ belonging to the Besov space $\dot{B}^{1,\infty}_{1}$, then $f\in BV$ and we have the inclusion $\dot{B}^{1,\infty}_{1}\subseteq BV$.
\end{Proposition}
\textit{\textbf{Proof.}}
Let $f\in \dot{B}^{1,\infty}_{1}(\mathbb{Z}_{2})$ and let us fix $|y|_{2}=2^{-m}$. We have to prove the following estimation for all $m>0$ 
\begin{equation*}
I=\int_{\mathbb{Z}_{2}}|f(x+y)-f(x)|dx\leq C\,2^{-m}.
\end{equation*}
Using the Littlewood-Paley decomposition given in (\ref{LPP}), we will work on the formula below
\begin{equation*}
I=\left\|\left(S_{0}f(x+y)+\sum_{j\geq 0}\Delta_{j}f(x+y)\right)-\left(S_{0}f(x)+\sum_{j\geq 0}\Delta_{j}f(x)\right)\right\|_{L^1}
\end{equation*}
Then, by the dyadic block's properties we have to study
\begin{equation}\label{Dosadico5}
I \leq \left\|S_{m}f(x+y)-S_{m}f(x)\right\|_{L^1}+\sum_{j=m+1}^{+\infty}\left\|\Delta_{j}f(x+y)-\Delta_{j}f(x)\right\|_{L^1}.
\end{equation}
We estimate this inequality with the two following lemmas.
\begin{Lemme}
The first term in (\ref{Dosadico5}) is identically zero. 
\end{Lemme}

\textit{\textbf{Proof.}}
Since we have fixed $|y|_{2}=2^{-m}$, then for $x\in Q_{m,k}$, we have $x+y\in Q_{m,k}$ with $k=\{0,...,2^{m}-1\}$. Applying the operators $S_{m}$ to the functions $f(x+y)$ and $f(x)$ we get the desired result.
\begin{flushright}{$\blacksquare$}\end{flushright} 

The second term in (\ref{Dosadico5}) is treated by the next lemma.
\begin{Lemme}
Under the hypothesis of proposition \ref{Dosadico} and for $|y|_{2}=2^{-m}$ we have
$$\sum_{j=m+1}^{+\infty}\left\|\Delta_{j}f(x+y)-\Delta_{j}f(x)\right\|_{L^1}\leq C\,2^{-m}.$$
\end{Lemme}
\textit{\textbf{Proof.}}
Indeed, 
$$\sum_{j=m+1}^{+\infty}\left\|\Delta_{j}f(x+y)-\Delta_{j}f(x)\right\|_{L^1}\leq 2 \sum_{j=m+1}^{+\infty}\left\|\Delta_{j}f\right\|_{L^1}.$$
We use now the fact $\|\Delta_{j}f\|_{L^1}\leq C\,2^{-j}$ for all $j$, since $f\in \dot{B}^{1,\infty}_{1}$, to get $$\sum_{j=m+1}^{+\infty}\left\|\Delta_{j}f(x+y)-\Delta_{j}f(x)\right\|_{L^1}\leq C\,2^{-m}.$$
\begin{flushright}{$\blacksquare$}\end{flushright} 
With these two lemmas, and getting back to (\ref{Dosadico5}), we deduce the following inequality for all $y\in \mathbb{Z}_{2}$:
$$\int_{\mathbb{Z}_{2}}|f(x+y)-f(x)|dx\leq C\,|y|_{2}$$
and this concludes the proof of proposition \ref{Dosadico}.
\begin{flushright} {$\blacksquare$}\end{flushright}

Our second step in order to prove theorem \ref{IdentityBVBesov} is the next result.
\begin{Proposition}\label{Unoadico}
In $\mathbb{Z}_{2}$ we have the inclusion $BV(\mathbb{Z}_{2})\subseteq \dot{B}^{1,\infty}_{1}(\mathbb{Z}_{2})$.
\end{Proposition}
\textit{\textbf{Proof.}}  Observe that we can characterize the Besov space $\dot{B}^{1,\infty}_{1}(\mathbb{Z}_{2})$ by the condition
$$\|f(\cdot+y)+f(\cdot-y)-2f(\cdot)\|_{L^1}\leq C\,|y|_{2}, \qquad \forall y\neq 0.$$
Let $f$  be a function in $BV(\mathbb{Z}_{2})$, then we have
$$\|f(\cdot+y)-f(\cdot)\|_{L^1}\leq C\;|y|_{2}.$$
Summing $\|f(\cdot-y)-f(\cdot)\|_{L^1}$ in both sides of the previous inequality we obtain
$$\|f(\cdot+y)-f(\cdot)\|_{L^1}+\|f(\cdot-y)-f(\cdot)\|_{L^1}\leq C\;|y|_{2}+\|f(\cdot-y)-f(\cdot)\|_{L^1}$$
and by the triangular inequality we have
$$\|f(\cdot+y)+f(\cdot-y)-2f(\cdot)\|_{L^1}\leq C\;|y|_{2}+\|f(\cdot-y)-f(\cdot)\|_{L^1}$$
We thus obtain
$$\|f(\cdot+y)+f(\cdot-y)-2f(\cdot)\|_{L^1}\leq 2C\;|y|_{2}.$$
\begin{flushright} {$\blacksquare$}\end{flushright}
We have proved, in the setting of the $2$-adic group $\mathbb{Z}_2$, the inequalities 
$$C_1\|f\|_{\dot{B}^{1,\infty}_1}\leq \|f\|_{BV} \leq C_2\|f\|_{\dot{B}^{1,\infty}_1},$$
so the theorem \ref{IdentityBVBesov} follows.
\begin{flushright}{$\blacksquare$}\end{flushright} 
\section{Improved Sobolev inequalities, $BV$ space and proof of theorem \ref{MainTheorem}}\label{ISPadique00}

We do not give here a global treatment of the family of inequalities of type (\ref{ISI2}); instead we focus on the next inequality
\begin{equation}\label{ISPadique0}
\|f\|_{L^2}^{2}\leq C\|f\|_{BV}\|f\|_{\dot{B}^{-1,\infty}_{\infty}}
\end{equation}
and we want to know if this estimation is true in a $2$-adic framework. Since in the $\mathbb{Z}_2$ setting we have the identification $\|f\|_{BV}\simeq\|f\|_{\dot{B}^{1,\infty}_{\infty}}$, the estimation (\ref{ISPadique0}) becomes
\begin{equation}\label{ISPadique1}
\|f\|_{L^2}^{2}\leq C\|f\|_{\dot{B}^{1,\infty}_{1}}\|f\|_{\dot{B}^{-1,\infty}_{\infty}}.
\end{equation}
This remark lead us to the theorem \ref{MainTheorem} which states that the previous inequalities are false.\\

\textit{\textbf{Proof.}} We will construct a counterexample by means of the Littlewood-Paley decomposition, so it is worth to recall very briefly the dyadic bloc characterization of the norms involved in inequality (\ref{ISPadique1}). For the $L^2$ norm we have $\|f\|_{L^2}^2=\sum_{j \in \mathbb{N}} \|\Delta_j f\|_{L^2}^2$, while for the Besov spaces $\dot{B}^{1,\infty}_{1}$ and $\dot{B}^{-1,\infty}_{\infty}$ we have
\begin{eqnarray*}
\|f\|_{\dot{B}^{1,\infty}_{1}}\; & = &\underset{ j\in \mathbb{N}}{\sup}\;2^{j}\|\Delta_j f\|_{L^1} \quad \mbox{ and}\\
\|f\|_{\dot{B}^{-1,\infty}_{\infty}}& = &\underset{ j\in \mathbb{N}}{\sup}\;2^{-j}\|\Delta_j f\|_{L^\infty}.
\end{eqnarray*}

We construct a function $f:\mathbb{Z}_2\longrightarrow \mathbb{R}$ by considering his values over the dyadic blocs and we will use for this the sets $Q_{j,k}$ defined in (\ref{QJKset}). First fix $\alpha$ and $\beta$ two non negative real numbers and $j_0, j_1$ two integers such that $0\leq j_0\leq j_1$ with the condition
$$2^{2j_0}\leq \frac{\beta}{\alpha}.$$
Now define $N_j$ as a function of $\alpha$ and $\beta$:
\begin{eqnarray}\label{ISPadique3}
N_j=2^j \quad\mbox{ if } 0\leq j \leq j_0 & \mbox{and} & N_{j}=\frac{\beta}{\alpha} 2^{-j}\leq 2^j \quad\mbox{ if } j_0< j \leq j_1.
\end{eqnarray}

and write
\begin{equation*}
\Delta_{j}f(x)=
\begin{cases}
\alpha 2^{j} &\text{over}\quad  Q_{j+1,0},\\
-\alpha 2^{j} &\text{over}\quad  Q_{j+1,1},\\
\alpha 2^{j} &\text{over}\quad  Q_{j+1,2},\\
-\alpha 2^{j} &\text{over}\quad  Q_{j+1,3},\\
&\vdots\\
\alpha 2^{j} &\text{over} \quad Q_{j+1,2N_j-2},\\
-\alpha 2^{j} &\text{over}\quad  Q_{j+1,2N_j-1},\\
0 &\text{elsewhere}.
\end{cases}
\end{equation*}

Once this function is fixed, we compute the following norms
\begin{itemize}
\item $\|\Delta_j f\|_{L^1}=\sum_{k=0}^{N_j}\alpha 2^{j}2^{-j}=\alpha N_j$,
\item $\|\Delta_j f\|_{L^\infty}=\alpha 2^{j}$,
\item $\|\Delta_j f\|_{L^2}^2=\sum_{k=0}^{N_j}\alpha^2 2^{2j}2^{-j}=\alpha^2 2^j N_j$,
\end{itemize}
and we build from these quantities the Besov and Lebesgue norms in the following manner:
\begin{enumerate}
\item[1)] For the Besov space $\dot{B}^{-1,\infty}_{\infty}$:\\

$\|f\|_{\dot{B}^{-1,\infty}_{\infty}}=\underset{0\leq j\leq j_1}{\sup}2^{-j}\alpha 2^j=\alpha$,
\item[2)] For the Besov space $\dot{B}^{1,\infty}_{1}$:\\

By the definition (\ref{ISPadique3}) of $N_j$ we have $2^j \|\Delta_j f\|_{L^1}=2^j \alpha N_j=2^{2j}\alpha$ if $0\leq j \leq j_0$ and $2^j \|\Delta_j f\|_{L^1}=\beta$ if $j_0<j\leq j_1$. Since $2^{2j_0}\leq \frac{\beta}{\alpha}$ we have:
$$\|f\|_{\dot{B}^{1,\infty}_{1}}=\beta.$$

\item[3)] For the Lebesgue space $L^2$:
\begin{eqnarray*}
\|f\|_{L^2}^2&=&\sum_{j=0}^{j_1} \alpha^2 2^j N_j=\sum_{j=0}^{j_0} \alpha^2 2^{2j} +\sum_{j>j_0}^{j_1} \alpha^2 2^j\frac{\beta}{\alpha} 2^{-j}= \sum_{j=0}^{j_0} \alpha^2 2^{2j}+ (j_1- j_0)\alpha \beta\\
&=& \alpha \beta \left( \frac{\alpha}{\beta}\sum_{j=0}^{j_0}  2^{2j}+ (j_1- j_0)\right). 
\end{eqnarray*}
With the condition $2^{2j_0}\leq \frac{\beta}{\alpha}$, we obtain from the previous formula that
$$\|f\|_{L^2}^2\simeq \alpha \beta (j_1- j_0)=\|f\|_{\dot{B}^{1,\infty}_{1}}\|f\|_{\dot{B}^{-1,\infty}_{\infty}}(j_1-j_0).$$
\end{enumerate}
Thus, getting back to (\ref{ISPadique1}) and therefore to (\ref{ISPadique0}), we have for an universal constant $C$ the inequality
\begin{eqnarray*}
\|f\|_{\dot{B}^{1,\infty}_{1}}\|f\|_{\dot{B}^{-1,\infty}_{\infty}}(j_1-j_0) & \leq & C \|f\|_{\dot{B}^{1,\infty}_{1}}\|f\|_{\dot{B}^{-1,\infty}_{\infty}}\\[5mm]
\iff (j_1-j_0) &\leq& C,
\end{eqnarray*}
which is false since we can freely choose the values of $j_1$ and $j_0$. The theorem \ref{MainTheorem} is proved.
\begin{flushright} {$\blacksquare$}\end{flushright}


\begin{flushright}
\begin{minipage}[r]{80mm}
Diego \textsc{Chamorro}\\[5mm]
Laboratoire d'Analyse et de Probabilités\\ 
Université d'Evry Val d'Essonne \& ENSIIE\\[2mm]
1 square de la résistance,\\
91025 Evry Cedex\\[2mm]
diego.chamorro@m4x.org
\end{minipage}
\end{flushright}

\end{document}